\documentclass[a4paper]{article}
\usepackage{amsmath,amssymb}
\usepackage[latin1]{inputenc}
\usepackage[dvips]{graphics}
\usepackage{mathrsfs}

\newcommand{\reg}{\operatorname{Reg}}
\newcommand{\gerst}{\operatorname{Gerst}}
\newcommand{\stlie}{\Sigma_t\operatorname{Lie}}
\newcommand{\com}{\operatorname{Com}}
\newcommand{\sym}{\mathfrak{S}}
\newcommand{\ac}{\mathsf{A}}
\newcommand{\anti}{\mathscr{A}}

\newtheorem{theorem}{Theorem}[section] 
\newtheorem{proposition}[theorem]{Proposition} 
\newtheorem{conjecture}[theorem]{Conjecture} 
 
\newtheorem{lemma}[theorem]{Lemme}

\newenvironment{proof}{\begin{trivlist}\item{\bf{Preuve.}}}
  {\hfill\rule{2mm}{2mm}\end{trivlist}}

\title{Sur le nombre de réflexions pleines dans les groupes de Coxeter
  finis}
\author{F. Chapoton}
\date{\today}

\begin{document}

\maketitle

\begin{abstract}
  On considère différents aspects d'une formule dans les groupes de
  Coxeter finis. 
\end{abstract}

\setcounter{section}{-1}

\section{Introduction}

Cet article tourne autour d'une formule qui définit, pour chaque
groupe de Coxeter fini $W$, un entier positif $f_W$ dépendant
seulement de la donnée des exposants de ce groupe de Coxeter. On peut
vérifier facilement, au cas par cas, que cet entier est le nombre de
réflexions dans $W$ qui sont pleines, \textit{i.e.} dont toutes les
décompositions réduites font intervenir tous les générateurs de
Coxeter.

Cet article comprend deux parties de nature différentes. Dans la
première partie, on cherche à obtenir une catégorification de la
formule, c'est à dire à l'interpréter comme une égalité de dimensions
provenant d'un isomorphisme entre deux modules sur le groupe de
Coxeter. On obtient une conjecture qui décrit précisément les modules
qui doivent entrer en jeu, puis on démontre cette conjecture dans les
cas des types $A$,$B$ et $I$.

La seconde partie est consacrée à une autre apparition de la formule
dans le contexte des systèmes de racines. Dans le cas des systèmes de
racines, les réflexions sont en bijection avec les racines positives.
Par cette bijection, les réflexions pleines correspondent aux racines
positives qui sont pleines, au sens où leur expression dans la base
des racines simples n'a pas de coefficient nul. Par une dualité
conjecturale sur l'ensemble des antichaînes du poset des racines
positives, les racines pleines devraient être en bijection avec les
antichaînes sans racines simples de cardinal maximal. On montre que le
nombre de telles antichaînes est bien égal au nombre de réflexions
pleines. On propose ensuite une conjecture reliant le polynôme $H$ qui
énumère les antichaînes selon leur cardinal et le nombre de racines
simples qu'elles contiennent à un polynôme $F$ similaire introduit
précédemment \cite{enum}. Par définition, un des coefficients de $H$
est donné par la formule qui nous intéresse ici.

\section{Réflexions pleines}

Soit $W$ un groupe de Coxeter fini de rang $n$ et
$S=\{s_1,\dots,s_n\}$ l'ensemble des réflexions simples de $W$. Une
réflexion $\sigma$ dans $W$ est dite \textit{pleine} si toute
décomposition réduite de $\sigma$ fait intervenir tous les éléments de
$S$. Dans le cas où $W$ est le groupe de Weyl d'un système de racines
cristallographique, les réflexions pleines correspondent aux racines
positives de support plein, \textit{i.e.} dont l'expression dans la
base des racines simples n'a pas de coefficient nul.

On vérifie aisément au cas par cas la proposition suivante, en
utilisant par exemple les tables de \cite{bourbaki}.

\begin{proposition}
  Le nombre $f_W$ de réflexions pleines dans $W$ est donné par la
  formule
  \begin{equation}
    \label{nombre_pleines}
    f_W=\frac{1}{|W|} (n h) \prod_{i=2}^{n} (e_i-1),
  \end{equation}
  où $|W|$ est l'ordre du groupe $W$, $h$ est le nombre de Coxeter et
  $e_1,\dots,e_n$ sont les exposants de $W$. Explicitement, on obtient :
  \begin{center}
    \begin{tabular}{|l|l|l|l|l|l|l|l|l|l|l|}
      \hline
      $A_n$ & $B_n$ & $D_n$ & $E_6$ & $E_7$ & $E_8$ &
      $F_4$ & $G_2$ & $H_3$ & $H_4$ & $I_2(h)$\\
      \hline
      $1$ & $n$ & $n-2$ & $7$ & $16$ & $44$ & $10$ & $4$ & $8$ & $42$ & $h-2$\\
      \hline
    \end{tabular}
  \end{center}
\end{proposition}

La formule (\ref{nombre_pleines}) peut se mettre sous la forme
suivante, plus suggestive.
\begin{equation}
  \label{ligne}
  (n h) \prod_{i=2}^{n} (e_i-1) = f_W |W|. 
\end{equation}

Il est alors naturel de chercher une interprétation de la formule
(\ref{ligne}) en termes d'un isomorphisme entre $W$-modules ou d'une
égalité entre caractères du groupe $W$. Comme $f_W$ est un entier, on
peut interpréter le membre de droite comme une somme directe de copies
de la représentation régulière $\reg$ de $W$. Le membre de gauche est
plus subtil.

On appelle \textit{racine} un demi-espace délimité par un des
hyperplans de $W$ dans l'espace euclidien. Soit $R$ la représentation
de $W$ sur l'ensemble des racines. Par une formule classique, la
dimension de $R$ est $n h$. Ceci fournit donc un $W$-module
susceptible de remplacer le premier facteur du membre de gauche de la
formule (\ref{ligne}).

Il reste donc à décrire un candidat pour le second facteur. Soit $G$
la re\-présentation de $W$ sur la cohomologie du complémentaire du
complexifié de l'arrangement d'hyperplans associé à $W$. Alors d'après
\cite{brieskorn}, la dimension graduée du $W$-module gradué $G$ est
donnée par
\begin{equation}
  \label{briesform}
  \sum_{k=0}^{n} \dim G_k (-t)^k = \prod_{i=1}^{n}(1-t e_i).
\end{equation}

On utilise alors \cite[Lemma 3.13]{orter} qui donne une différentielle
acyclique naturelle $\partial$ sur l'algèbre de Orlik-Solomon d'un
arrangement d'hyperplans non vide \cite{orsol}. Le Lemme suivant en
est une conséquence immédiate.
\begin{lemma}
  \label{crux}
  Le caractère de $G$ est divisible par $1-t$ et le quotient est le
  caractère d'une représentation graduée.
\end{lemma}

Remarque : Comme $\partial$ est en fait une dérivation pour le produit
en cohomologie, le quotient $\frac{G}{1-t}$ a une structure d'algèbre
graduée, liée dans le cas des groupes symétriques à la cohomologie des
espaces de modules de courbes de genre $0$ (voir \cite{getzler}).

Soit $G'$ le $W$-module virtuel obtenu en faisant $t=1$ dans le
$W$-module gradué $\frac{G}{1-t}$. Par la formule (\ref{briesform}), la
dimension (virtuelle) de $G'$ est $\prod_{i=2}^{n}(1-e_i)$.

On a donc défini des $W$-modules $R$ et $\reg$ et un $W$-module
virtuel $G'$ qui vérifient une égalité de dimensions équivalente à la
formule (\ref{ligne}) :
\begin{equation}
   (\dim R)(\dim G')= (-1)^{n-1} f_W \dim(\reg).
\end{equation}

Cette égalité devrait être une conséquence de la conjecture suivante.
\begin{conjecture}
  \label{principale}
  On a une égalité de caractères :
  \begin{equation}
    R \otimes G' = (-1)^{n-1} f_W \reg.
  \end{equation}
\end{conjecture}

On peut peut-être espérer un énoncé plus précis, comme l'existence
d'une différentielle sur le $W$-module gradué
\begin{equation}
  R \otimes \frac{G}{1-t},
\end{equation}
dont l'homologie serait concentrée en degré $n-1$ et isomorphe en ce
degré à une somme de $f_W$ copies de la représentation régulière.

Si une telle différentielle existe, il doit exister un $W$-module dont
le caractère est donné par 
\begin{equation}
  \frac{1}{1-t}\left(R \otimes \frac{G}{1-t}-(-t)^{n-1} f_W \reg\right).
\end{equation}
On peut vérifier que ceci est vrai pour les groupes symétriques de
petit rang.

Les trois sections suivantes sont consacrées à la preuve de la
conjecture \ref{principale} pour les groupes symétriques, les groupes
hyperoctaédraux et les groupes diédraux respectivement.

\section{Le cas des groupes symétriques}

Dans cette section, on démontre la conjecture \ref{principale} dans le
cas du groupe de Coxeter de type $A_{n-1}$, qui est le groupe
symétrique $\sym_n$ sur $n$ lettres. Une autre preuve est probablement
possible, dans l'esprit de celle donnée plus loin pour le type $B$.

\subsection{Le caractère de $R$}

Étudions d'abord le caractère de $R$. Soit $N$ la représentation
naturelle de dimension $n$, \textit{i.e.} l'action de $\sym_n$ par
permutations de l'ensemble $\{1,\dots,n\}$. Alors on a clairement un
isomorphisme $N\otimes N \simeq N \oplus R$. Le caractère $\chi_N$ de
$N$ est facile à décrire. Si $\lambda$ est une partition de $n$, on
note $m_\lambda$ le nombre de parts de taille $1$ dans $\lambda$. Sur
la classe de conjugaison $C_\lambda$ associée à une partition
$\lambda$, on a $\chi_N(C_\lambda)=m_\lambda$. Par conséquent, on
obtient $\chi_R(C_\lambda)=m_\lambda^2-m_\lambda$.

\begin{lemma}
  Le caractère $\chi_R$ s'annule sur la classe de conjugaison
  $C_\lambda$ si et seulement si $m_\lambda$ est au plus égal à $1$.
\end{lemma}

\subsection{Le caractère de $G'$}

Pour démontrer que le caractère de $R \otimes G'$ est un multiple du
caractère de la représentation régulière, il suffit donc de montrer
que le caractère $\chi_{G'}$ de $G'$ vérifie la condition suivante :
\begin{equation}
  \text{Si }m_\lambda \geq 2 
  \text{ et }\lambda \not= 1^n 
  \text{ alors } \chi_{G'}(C_\lambda)=0.
\end{equation}
Par la définition de $G'$, ceci est équivalent à la condition
\begin{equation}
  \text{Si }m_\lambda \geq 2 
  \text{ et }\lambda \not= 1^n 
  \text{ alors } \chi_{G}(C_\lambda)\text{ a une racine double en }t=1.
\end{equation}

On passe au langage des fonctions symétriques, en identifiant un
caractère à une fonction symétrique de la manière habituelle. Pour
$n\geq 1$, on note $p_n$ la fonction symétrique ``somme des
puissances'' d'ordre $n$. L'assertion précédente est donc équivalente
au Lemme suivant.
\begin{lemma}
  \label{bonzero}
  La valeur de
  \begin{equation}
    \frac{1}{1-t}\partial^2_{p_1} \chi_{G}
  \end{equation}
  en $t=1$ est proportionnelle à la fonction symétrique $p_1^n$.
\end{lemma}

La preuve de ce Lemme est obtenue dans la section suivante.

\subsection{Séries génératrices}

Soit $\gerst$ la série génératrice des caractères $\chi_G$ :
\begin{equation}
  \gerst=\sum_{n\geq 1}\chi_{G^n},
\end{equation}
où $G^n$ est le module $G$ pour le groupe $\sym_n$. Par convention,
$G^1$ est le module trivial.

On dispose d'une description de $\gerst$ par le biais de la théorie
des opérades. En effet, le complémentaire du complexifié de
l'arrangement de type $A_{n-1}$ est homotope à l'espace des petits
disques, formé par l'ensemble des plongements disjoints de $n$ disques
dans le disque unité du plan complexe. Il en résulte un isomorphisme
en homologie. Mais les petits disques ont une structure d'opérade
topologique et leur homologie est l'opérade dite de Gerstenhaber, voir
\cite{voronov}. On sait par ailleurs que le foncteur analytique
sous-jacent à l'opérade de Gerstenhaber est le composé des foncteurs
analytiques $\com$ sous-jacent à l'opérade des algèbres commutatives
et du foncteur $\stlie$ sous-jacent à la suspension de l'opérade des
algèbres de Lie, voir par exemple \cite{markl}.

Soit donc $\com$ la fonction symétrique $\sum_{n\geq 1} h_n$ où $h_n$
est la fonction symétri\-que complète, et soit $\stlie$ la fonction symétrique
\begin{equation}
  \sum_{n \geq 1}\frac{(-t)^{n-1}}{n}\sum_{d|n} \mu(d) p_d^{n/d},
\end{equation}
où $\mu$ est la fonction de Möbius.

La traduction en termes de fonctions symétriques de la relation de
composition entre foncteurs analytiques décrite ci-dessus est l'énoncé
suivant.

\begin{proposition}
  On a la relation pléthystique :
  \begin{equation}
    \gerst=\com \circ \left(\stlie\right).
  \end{equation}
\end{proposition}

De plus, on calcule facilement 
\begin{equation}
  \partial_{p_1}\com =1+\com \text{ et }
  \partial_{p_1}\stlie =\frac{1}{1+p_1 t}.
\end{equation}

On en déduit, en utilisant le fait que $\partial_{p_1}$ est une
dérivation pour le pléthysme, que
\begin{equation}
   \partial^2_{p_1} \gerst =
   (1-t) \frac{1}{(1+p_1 t)^2}(1+\com) \circ \left(\stlie\right).
\end{equation}

Par le Lemme \ref{crux}, la valeur de $\com \circ
\left(\stlie\right)$ en $t=1$ est $p_1$. On a donc
\begin{equation}
   \left(\frac{1}{1-t}\partial^2_{p_1} \gerst \right)\bigg{|}_{t=1}=
   \frac{1}{(1+p_1)^2}(1+p_1).
\end{equation}

Comme cette expression est une fonction de $p_1$ seulement, on a
démontré le Lemme \ref{bonzero} et donc la conjecture \ref{principale}
pour les groupes de Coxeter de type $A$.
\section{Groupes hyperoctaédraux}

On considère ici la cas du groupe de Coxeter $W$ de type $B_n$. On
renvoie à \cite{lehrer} pour la description des classes de conjugaison
de $W$ par des paires de partitions. 

Soit $R$ le module des racines. Si $(\epsilon_i)_{i \in\{1,\dots,n\}}$
est une base orthonormale, les racines sont $\pm\epsilon_i$ et $\pm
\epsilon_i \pm \epsilon_j$ pour $i\not=j$. Commençons par décrire les
classes de conjugaison non triviales de $W$ sur lesquelles le
caractère de $R$ est non nul, \textit{i.e.} qui ont au moins un point
fixe dans $R$. Si un élément $g$ du groupe fixe la racine
$\epsilon_i$, il doit avoir un point fixe $i$, \textit{i.e.} un cycle
positif de longueur $1$. Si $g$ fixe une racine
$\epsilon_i-\epsilon_j$, alors il doit vérifier $g(i)=-j$ et
$g(j)=-i$. De même, si $g$ fixe une racine $\epsilon_i+\epsilon_j$, il
doit vérifier $g(i)=j$ et $g(j)=i$. Dans ces deux cas, $g$ a un cycle
positif de longueur $2$.

En conclusion, on a obtenu le Lemme suivant.
\begin{lemma}
  Le caractère de $g$ sur $R$ est non nul si et seulement si $g$ a un
  cycle positif de longueur $1$ ou un cycle positif de longueur $2$.
\end{lemma}

On utilise ensuite la description de la cohomologie due à Lehrer
\cite[Thm. 5.6]{lehrer}. On change $t$ en $-t$ dans les résultats de
Lehrer pour respecter nos conventions. En particulier, la valeur du
caractère de la cohomologie sur un élément $g$,
\begin{equation}
  \chi(g)=\sum_{k=0}^{n} tr_g(G_k) (-t)^k,
\end{equation}
est décrite dans \cite{lehrer} comme un produit de facteurs explicites
ne dépendant que de la décomposition en cycles positifs et négatifs de
$g$.

Si un élément $g$ de $W$ ayant un cycle positif de longueur $1$ est
non trivial, alors $g$ possède soit un cycle négatif soit un cycle
positif de longueur au moins $2$, \textit{i.e.} $g$ possède au moins
deux cycles de longueurs ou signes différents. D'après les résultats
de \cite{lehrer}, le polynôme $\chi(g)$ est divisible par $1-t$ autant
de fois qu'il y a de cycles dans $g$. On en déduit que le caractère de
$G'$ est nul sur la classe de $g$.

Si maintenant $g$ a un cycle positif de longueur $2$, on déduit de
\cite{lehrer} que le caractère de $G$ sur la classe de $g$ est
divisible par $(1-t)^2$. Par conséquent, le caractère de $G'$ sur la
classe de $g$ est nul.

On a donc montré le Lemme suivant.
\begin{lemma}
  Soit $g$ un élément non trivial de $W$ qui a un cycle positif de
  longueur $1$ ou un cycle positif de longueur $2$. Alors le caractère
  de $G'$ est nul sur la classe de conjugaison de $g$.
\end{lemma}

Les deux Lemmes précédents entraînent la conjecture \ref{principale}
pour le type $B_n$.

\section{Groupes diédraux}

On considère maintenant le cas des groupes diédraux. Soit $W$ un
groupe de type $I_2(h)$.

Si $h$ est impair, alors le module $R$ donné par l'action de $W$ sur
les racines est isomorphe à la représentation régulière $\reg$, donc
la conjecture \ref{principale} est trivialement vraie, car les
dimensions sont égales.

Supposons donc $h$ pair. Dans ce cas, un élément non trivial de $W$ a
un trace non nulle dans le module $R$ des racines si et seulement si
c'est une réflexion. Il faut donc montrer que le caractère de $G'$
s'annule sur les réflexions. Soit donc $H_0$ un hyperplan fixé parmi
les hyperplans de $W$ et soient $H_{1},\dots,H_{h-1}$ les autres
hyperplans dans un ordre cyclique. Soit $\sigma$ la réflexion par
rapport à $H_0$.

Comme le rang du groupe de Coxeter $W$ est $2$, on dispose d'une
description explicite du module gradué $\frac{G}{1-t}$ comme suit.

En degré $0$, le module a pour base la fonction $1$, donc la trace de
$\sigma$ est $1$.

En degré $1$, le module est le noyau de la différentielle $\partial$
de degré $-1$ sur la cohomologie en degré $1$. On rappelle
\cite{orter} que cette différentielle est définie par
$\partial(\omega_H)=1$ pour tout hyperplan $H$ de l'arrangement, où
$\omega_H$ est la forme différentielle logarithmique associée à $H$.
Par conséquent, ce module a une base $O_i=\omega_{H_0}-\omega_{H_i}$
pour $i=1,\dots,h-1$.  L'action de $\sigma$ est donnée par
$\sigma(O_i)=O_{h-i}$. Il y a un seul point fixe qui est $O_{h/2}$,
donc la réflexion $\sigma$ a pour trace $1$.

Au total, la trace de la réflexion $\sigma$ sur $G'$ est donc nulle,
ce qui démontre la conjecture \ref{principale} pour les groupes
diédraux.

\section{Antichaînes sans racines simples}

On considère maintenant une autre apparition de la formule
(\ref{nombre_pleines}) dans le contexte des systèmes de racines.

Soit $\Phi$ un système de racines fini de rang $n$, $\Phi_{\geq 0}$
l'ensemble des racines positives et $\Pi=\{\alpha_i\}_{i \in I}$
l'ensemble des racines simples. Soit $W$ le groupe de Weyl associé. Il
existe un ordre partiel naturel $\leq$ sur $\Phi_{\geq 0}$ défini par
$\alpha \leq \beta$ si $\beta-\alpha$ est une combinaison linéaire à
coefficients positifs de racines simples.

Une \textit{antichaîne} dans un poset $P$ est une partie de $P$ formée
d'éléments tous incomparables pour $\leq$. On note $\anti(P)$
l'ensemble des antichaînes de $P$. On considère implicitement par la
suite, sauf précision contraire, que le poset considéré est
$\Phi_{\geq 0}$ pour la relation $\leq$.

Les antichaînes dans ce poset des racines positives ont été beaucoup
étudiés récemment \cite{athaBull,athaTAMS,ysystem}. En particulier, on
sait que leur nombre est égal au nombre de Catalan généralisé associé
au système de racine $\Phi$. Le lemme suivant est bien connu.
\begin{lemma}
  \label{nantichaine}
  Le cardinal maximal d'une antichaîne est $n$. L'unique telle
  antichaîne est $\ac=\{\alpha_i\}_{i\in I}$.
\end{lemma}

On constate au cas par cas (voir \cite{athaTAMS}) que le nombre
d'antichaînes de cardinal $k$ est égal au nombre d'antichaînes de
cardinal $n-k$, pour tout $k$. Conjecturalement, ceci doit être une
conséquence de l'existence d'une dualité sur l'ensemble des
antichaînes qui envoie les antichaînes de cardinal $k$ sur les
antichaînes de cardinal $n-k$. Une dualité ayant cette propriété a été
construite par Panyushev pour les types $A$,$B$ et $C$ dans
\cite{pany}.

Soit $\ac$ une antichaîne. On définit le \textit{type} de $\ac$ noté
$t(\ac)$ qui est une sous-partie de l'ensemble des arêtes du diagramme
de Dynkin. Une arête $(i,j)$ est dans $t(\ac)$ si et seulement si il
existe une racine $\alpha$ dans $\ac$ telle que $i$ et $j$ soient dans
le support de $\alpha$.

Panyushev propose de chercher une dualité ayant des propriétés plus
fortes, \textit{cf.} sa conjecture \cite[Conj. 6.1]{pany}. On peut
reformuler en partie ses hypothèses sous la forme suivante : la
dualité doit préserver le type.

Si le type $t(\ac)$ est le diagramme de Dynkin complet, on dit que
$\ac$ est de type plein. Comme cas particulier de la conjecture
d'existence d'une dualité respectant le type, on a la conjecture
ci-dessous.

\begin{conjecture}
  \label{dual_racine_pleine}
  Il existe une bijection naturelle entre les antichaînes de cardinal
  $1$ de type plein et les antichaînes de cardinal $n-1$ de type
  plein.
\end{conjecture}

Une antichaîne $\ac$ de cardinal $1$ est juste une racine positive
$\alpha$. Une racine positive $\alpha$ est de type plein comme
antichaîne si et seulement si elle est de support plein, \textit{i.e.}
correspond à une réflexion pleine. Il y a donc une bijection entre
antichaînes de cardinal $1$ de type plein et réflexions pleines.

\begin{lemma}
  Une antichaîne $\ac$ de cardinal $n-1$ est de type plein si et
  seulement si elle ne contient pas de racine simple.
\end{lemma}
\begin{proof}
  C'est clairement vrai si $n=1$. Supposons donc $n\geq 2$. Si une
  antichaîne $\ac$ contient une racine simple $\alpha_i$, elle ne peut
  contenir aucune autre racine ayant $\alpha_i$ dans son support. Par
  conséquent toute arête de la forme $(i,j)$ ne peut pas être couverte
  par $\ac$. Une telle arête existe car $n$ est au moins $2$, donc
  $\ac$ n'est pas de type plein.
  
  Réciproquement, si $\ac$ est une antichaîne de cardinal $n-1$ de
  type non plein, il existe une arête non couverte par $\ac$.
  Considérons le diagramme de Dynkin non connexe obtenu en enlevant
  cette arête. Il existe une partition de ce diagramme en deux parties
  non vides de cardinal $p$ et $q$ avec $p+q=n$, sans arêtes entre
  elles. Alors $\ac$ est l'union disjointe de deux antichaînes $\ac_1$
  et $\ac_2$ contenues respectivement dans les posets associés à ces
  deux diagrammes. Par le Lemme \ref{nantichaine}, le cardinal de
  $\ac_1$ est inférieur ou égal à $p$ et celui de $\ac_1$ est
  inférieur ou égal à $q$. Quitte à échanger les indices $1$ et $2$,
  on peut supposer que $\ac_1$ est une antichaîne de cardinal $p$. Par
  le Lemme \ref{nantichaine}, $\ac_1$ est donc formée de racines
  simples. Donc $\ac$ contient au moins une racine simple.
\end{proof}

La conjecture \ref{dual_racine_pleine} se reformule donc ainsi.
\begin{conjecture}
  Il existe une bijection naturelle entre les racines pleines et les
  antichaînes de cardinal $n-1$ sans racines simples.
\end{conjecture}

Vérifions que cette relation est vraie au niveau des cardinaux de ces
ensembles.

\begin{proposition}
  \label{antichaine_formule}
  Le nombre d'antichaînes de cardinal $n-1$ sans racines simples est
  égal au nombre de réflexions pleines, donc donné par la formule
  (\ref{nombre_pleines}).
\end{proposition}

\begin{proof}
  On utilise les polynômes de Narayana généralisés définis par
  \begin{equation}
    N_{\Phi}(x)=\sum_{k=0}^{n} n_k(\Phi) x^k,
  \end{equation}
  où $n_k(\Phi)$ est le nombre d'antichaînes de cardinal $k$ dans le
  poset $(\Phi_{\geq 0},\leq)$. On constate au cas par cas sur les
  résultats d'Athanasiadis \cite{athaTAMS} que ces polynômes sont
  symétriques, \textit{i.e.} vérifient la relation
  \begin{equation}
    N_{\Phi}(x)=x^n N_{\Phi}(1/x).
  \end{equation}

  Considérons maintenant les polynômes 
  \begin{equation}
    P_{\Phi}(x)=\sum_{k=0}^{n} p_k(\Phi) x^k,
  \end{equation}
  où $p_k(\Phi)$ est le nombre d'antichaînes pleines de cardinal $k$
  dans le poset $(\Phi_{\geq 0},\leq)$.  
  
  Les polynômes $N$ et $P$ sont en fait naturellement définis pour des
  diagrammes de Dynkin non nécessairement connexes et sont alors
  donnés par le produit des polynômes associés aux composantes
  connexes.

  Par définition de la notion de type, on a la relation suivante entre
  les polynômes $N$ et $P$ :
  \begin{equation}
    N_{\Phi}=\sum_{E} P_{\Phi[E]},
  \end{equation}
  où la somme porte sur l'ensemble des parties $E$ de l'ensemble des
  arêtes du diagramme de Dynkin de $\Phi$ et $\Phi[E]$ désigne le
  diagramme de Dynkin non nécessairement connexe obtenu en ne gardant
  que les arêtes dans $E$. 
  
  Par conséquent, par inversion de Möbius, on peut exprimer les
  polynômes $P$ en fonction des polynômes $N$. On en déduit que
  \begin{equation}
    P_{\Phi}(x)=x^n P_{\Phi}(1/x).
  \end{equation}
  En particulier, le nombre des antichaînes de type plein de
  cardinal $n-1$ est égal au nombre de réflexions pleines.
\end{proof}
 
\section{Une conjecture énumérative}

On conjecture ici une relation entre le polynôme énumérateur $F$ des
associaèdres généralisés introduit dans \cite{enum} et un polynôme
énumérateur $H$ des antichaînes selon deux paramètres. Un des
coefficients de $H$ est le nombre de la formule
(\ref{nombre_pleines}).

Soit $\Phi$ un système de racines. Considérons le polynôme suivant :
\begin{equation}
  H = \sum_{k,\ell} h_{k,\ell}\, x^k y^\ell,
\end{equation}
où $h_{k,\ell}$ est le nombre d'antichaînes de cardinal $k$ contenant
$\ell$ racines simples. En particulier $h_{n-1,0}$ est le nombre
d'antichaînes de cardinal $n-1$ sans racines simples considéré
précédemment.

On rappelle brièvement la définition du polynôme $F$. Dans
\cite{ysystem}, Fomin et Zelevinsky ont introduit pour chaque système
de racines un complexe simplicial $\Delta(\Phi)$ dont l'ensemble des
sommets est l'ensemble des racines presque positives $\Phi_{\geq
  -1}=\Phi_{\geq 0} \sqcup (-\Pi)$. On définit le polynôme $F$ par la
formule
\begin{equation}
  F(x,y)=\sum_{k=0}^{n}\sum_{\ell=0}^{n}f_{k,\ell}x^k y^\ell,
\end{equation}
où $f_{k,\ell}$ est le nombre de simplexes de $\Delta(\Phi)$ ayant
exactement $k$ sommets dans $\Phi_{\geq 0}$ et $\ell$ sommets dans
$-\Pi$. Ce polynôme a été défini et étudié dans \cite{enum}.

\begin{conjecture}
  On a la relation suivante :
  \begin{equation}
    H(x,y)=(1-x)^n\, F(x/(1-x),xy/(1-x)).
  \end{equation}
\end{conjecture}

On peut aisément vérifier cette conjecture pour les systèmes de
racines de petit rang.

\bibliographystyle{plain}
\bibliography{greg}

\begin{thebibliography}{10}

\bibitem{athaBull}
C.~A. Athanasiadis.
\newblock Generalized {C}atalan numbers, {W}eyl groups and arrangements of
  hyperplanes,.
\newblock {\em Bull. London Math. Soc.}, 2004.

\bibitem{athaTAMS}
C.~A. Athanasiadis.
\newblock On a refinement of the generalized {C}atalan numbers for {W}eyl
  groups.
\newblock {\em Trans. A. M.S.}, 2004.

\bibitem{bourbaki}
N.~Bourbaki.
\newblock {\em \'{E}l\'ements de math\'ematique. {F}asc. {XXXIV}. {G}roupes et
  alg\`ebres de {L}ie. {C}hapitre {IV}, {C}hapitre {V}, {C}hapitre {VI}}.
\newblock Hermann, 1968.

\bibitem{brieskorn}
E.~Brieskorn.
\newblock Sur les groupes de tresses [d'apr\`es {V}. {I}. {A}rnold].
\newblock In {\em S\'eminaire Bourbaki, 24\`eme ann\'ee (1971/1972), Exp. No.
  401}, pages 21--44. Lecture Notes in Math., Vol. 317. Springer, 1973.

\bibitem{enum}
F.~Chapoton.
\newblock Enumerative properties of generalized associahedra.
\newblock {\em Séminaire Lotharingien de Combinatoire}, 51, 2004.

\bibitem{ysystem}
S.~Fomin and A.~Zelevinsky.
\newblock {$Y$}-systems and generalized associahedra.
\newblock {\em Ann. of Math. (2)}, 158(3):977--1018, 2003.

\bibitem{getzler}
E.~Getzler.
\newblock Operads and moduli spaces of genus {$0$} {R}iemann surfaces.
\newblock In {\em The moduli space of curves (Texel Island, 1994)}, volume 129
  of {\em Progr. Math.}, pages 199--230. Birkh\"auser Boston, 1995.

\bibitem{lehrer}
G.~I. Lehrer.
\newblock On hyperoctahedral hyperplane complements.
\newblock volume~47 of {\em Proc. Sympos. Pure Math.}, pages 219--234. Amer.
  Math. Soc., 1987.

\bibitem{markl}
M.~Markl.
\newblock Distributive laws and {K}oszulness.
\newblock {\em Ann. Inst. Fourier (Grenoble)}, 46(2):307--323, 1996.

\bibitem{orsol}
P.~Orlik and L.~Solomon.
\newblock Combinatorics and topology of complements of hyperplanes.
\newblock {\em Invent. Math.}, 56(2):167--189, 1980.

\bibitem{orter}
P.~Orlik and H.~Terao.
\newblock {\em Arrangements of hyperplanes}, volume 300 of {\em Grundlehren der
  Mathematischen Wissenschaften}.
\newblock Springer-Verlag, 1992.

\bibitem{pany}
D.~I. Panyushev.
\newblock Ad-nilpotent ideals of a {B}orel subalgebra: generators and duality.
\newblock {\em J. Algebra}, 2004.

\bibitem{voronov}
A.~A. Voronov.
\newblock The {S}wiss-cheese operad.
\newblock In {\em Homotopy invariant algebraic structures (Baltimore, MD,
  1998)}, volume 239 of {\em Contemp. Math.}, pages 365--373. Amer. Math. Soc.,
  1999.

\end{thebibliography}

\end{document}